\documentclass[a4paper,10pt]{amsart}
\usepackage[ansinew]{inputenc}
\usepackage[all]{xy}
\usepackage{amsmath,latexsym,amssymb,verbatim,amsthm}

\input arrow.tex

\newcommand{\N}{\mathbb{N}}
\newcommand{\Z}{\mathbb{Z}}

\newtheorem{lema}{Lemma}[section]
\newtheorem{teo}[lema]{Theorem}
\newtheorem{pro}[lema]{Proposition}
\newtheorem{defi}[lema]{Definition}
\newtheorem{cor}[lema]{Corollary}
\newtheorem{com}[lema]{Remark}
\newtheorem{eje}[lema]{Example}
\newtheorem{den}[lema]{Notation}
\newtheorem{teo*}{Theorem}

\begin{document}

\title{Functorial Cartier Duality}

\author{Amelia \'{A}lvarez}
\address{Departamento de Matem\'{a}ticas, Universidad de Extremadura,
Avenida de Elvas s/n, 06071 Badajoz, Spain}
\email{aalarma@unex.es}

\author{Carlos Sancho}
\address{Departamento de Matem\'{a}ticas, Universidad de Salamanca,
Plaza de la Merced 1-4, 37008 Salamanca, Spain}
\email{mplu@usal.es}

\author{Pedro Sancho}
\address{Departamento de Matem\'{a}ticas, Universidad de Extremadura,
Avenida de Elvas s/n, 06071 Badajoz, Spain} \email{sancho@unex.es}

\subjclass[2000]{Primary 14L15. Secondary 20M07, 22D35}

\date{June, 2007}

\maketitle

\section*{Introduction}

In this paper we obtain the Cartier duality for $k$-schemes of 
commutative monoids functorially without providing the vector
spaces of functions with a topology (as in \cite[Expos\'{e} VII$_B$
by P. Gabriel, 2.2.1]{demazure}), generalizing a result for finite
commutative algebraic groups by M. Demazure $\&$ P. Gabriel
(\cite[II, $\S$1, 2.10]{Demazure}).

All functors we consider are functors defined over the category of
commutative $k$-algebras. Given a $k$-module $E$, we denote by
${\bf E}$ the functor of $k$-modules ${\bf E}(B):= E \otimes_k B$.
Given two functors of $k$-modules $F$ and $H$, ${\bf Hom}_k(F,H)$
will denote the functor of $k$-modules $${\bf Hom}_k(F,H)(B)={\rm
Hom}_{B}(F_{|B}, H_{|B})$$ where $F_{|B}$ is the functor $F$
restricted to the category of commutative $B$-algebras. We denote
$F^* := {\bf Hom}_k (F ,{\bf k})$. It holds that ${\bf E}^{**} =
{\bf E}$ (\cite[1.10]{Amel}). Given ${\bf E}$ there exists a
$k$-module $V$  such that ${\bf E} = {\bf V}^*$ if and only if $E$
is a projective $k$-module of finite type.

Given a functor of commutative $k$-algebras ${\mathcal A}$, we
define ${\rm Spec}\, {\mathcal A}$ to be the functor ${\rm Spec}\,
{\mathcal A}(B) := {\rm Hom}_{k-alg}({\mathcal A}, {\bf B})$. Let
$X={\rm Spec}\, A$ be a $k$-scheme and let $X^\cdot$ be its
functor of points ($X^\cdot(B):={\rm Hom}_{k-alg}(A,B)$). It is
easy to see that $X^\cdot ={\rm Spec}\, {\bf A}$. If ${\bf C}^*$
is a functor of commutative algebras, we prove (Theorem \ref{1.5})
that ${\rm Spec}\, {\bf C} ^*$ is a direct limit (of functor of
points) of finite $k$-schemes.

If ${\mathcal A} = {\bf A}$, in \cite[3.4]{Amel} we proved that
the linear enveloping of ${\rm Spec}\, {\mathcal A}$ is ${\mathcal
A}^*$, that is, $${\rm Hom}_{functors} ({\rm Spec}\, {\mathcal A},
H) = {\rm Hom}_k ({\mathcal A}^*,H)$$ for all dual functors of
$k$-modules $H$ (i.e., $H= F^*$). If in addition $G= {\rm Spec}\,
{\mathcal A}$ is an affine $k$-group, in \cite[5.4]{Amel} we proved that the algebraic enveloping of ${\rm Spec}\, {\mathcal A}$ is ${\mathcal
A}^*$ again, that is, $${\bf Hom}_{monoids} ({\rm Spec}\,
{\mathcal A}, {\mathcal B}) = {\bf Hom}_{k-alg} ({\mathcal A}^*,
{\mathcal B}) $$ for all dual functors of $k$-algebras ${\mathcal
B}$. As a consequence, we obtained that the category of $G$-modules is equivalent to the category of ${\mathcal A}^*$-modules.

In this paper we prove simultaneously the same results for
${\mathcal A} = {\bf A}, {\bf C}^*$ (Corollary \ref{2.3}, Proposition \ref{2.5}, Theorem \ref{Gmodulos}). In particular, $${\bf Hom}_{functors} ({\rm Spec}\, {\mathcal A}, {\bf k}) = {\bf Hom}_k ({\mathcal A}^*, {\bf k} ) = {\mathcal A}.$$ Moreover, if ${\rm Spec}\, {\mathcal A}$ is a functor of abelian monoids, then we get $$ {\bf Hom}_{monoids} ({\rm Spec}\,{\mathcal A}, {\bf k}) = {\bf Hom}_{k-alg} ({\mathcal A}^*,{\bf k}) = {\rm Spec}\, {\mathcal A}^* $$ where ${\bf k}$ is regarded as a functor of monoids with the multiplication
operation.

Given a functor of monoids $H$, we say that $H^* = {\bf
Hom}_{monoids} (H,{\bf k})$ is the dual functor of monoids of $H$
and we prove the following theorem.

\begin{teo*}
The category of abelian affine monoids $G={\rm Spec}\, A$ is
anti-equivalent to the category of functors ${\rm Spec}\, {\bf
C}^*$ of abelian monoids. The functors giving the anti-equivalence
are the ones that assign to each functor of monoids its dual functor of monoids.
\end{teo*}

In particular, $G^{**}=G$.

Finally, we get three immediate applications of the functorial
Cartier duality:
\begin{enumerate}
\item[(1)] The tangent space to $G^*$ at the identity element is
isomorphic to the vector space $V$ of linear functions of $G$.

\item[(2)] If ${\rm car}\, k =0$ and $G$ is a unipotent abelian
$k$-group, then it is isomorphic to the dual of $V$.

\item[(3)] If $G$ is semisimple, then the transposed Fourier transform of
$G^*$ is the inverse of the Fourier transform of $G$.
\end{enumerate}

\section{Spectrum of a functor of algebras}

The basic references for reading this paper are \cite{Amel} and
\cite{carlos}.

Let $k$ be a commutative ring with unit. If $E$ is a $k$-module we
will say that ${\bf E}$ is a cuasi-coherent $k$-module. The
category of $k$-modules is equivalent to the category of
cuasi-coherent $k$-modules (\cite[1.12]{Amel}). In particular,
${\rm Hom}_k ({\bf E},{\bf E'}) = {\rm Hom}_k (E,E')$. ${\bf E}^*
= {\bf Hom}_k ({\bf E},{\bf k})$ is the functor of points of the
scheme ${\rm Spec}\,S_k^\cdot E$ and  we say that it is a
$k$-module scheme. As it holds that ${\bf E}^{**} = {\bf E}$, the category of cuasi-coherent modules is anti-equivalent to the category of
module schemes (\cite[1.12]{Amel}) .

\begin{den}
Throughout this paper we assume that ${\bf A}^*$ is a commutative
$k$-algebra scheme (equivalently, $A$ is a cocommutative
$k$-coalgebra, \cite[4.2]{Amel}). From Proposition $1.5$ on, we
will always assume that $A$ is a projective $k$-module.
\end{den}

\begin{defi}
Given a functor of $k$-algebras $\mathcal{A}$, the functor
${\rm Spec}\, \mathcal{A}$, ``spectrum of $\mathcal{A}$'', is defined
to be $$({\rm Spec}\, \mathcal{A})(B) := {\rm Hom}_{k-alg}
(\mathcal{A}, {\bf B})$$ for each commutative $k$-algebra $B$.
\end{defi}

\begin{eje}
If $C$ is a commutative $k$-algebra, then $${\rm Spec}\, {\bf C} =
({\rm Spec}\, C)^\cdot .$$
\end{eje}

\begin{eje}
Let $X$ be a set. Let us consider the discrete topology on $X$.
Let $\mathcal{X}$ be the functor, which we will call the constant
functor $X$, defined by $${\mathcal X}(B) := {\rm Aplic}_{cont.}
({\rm Spec}\, B, X)$$ for each commutative $k$-algebra $B$.

Let ${\mathcal A}_X$ be the functor of algebras defined by
$${\mathcal A}_X(B) := {\rm Aplic}(X,B) = \prod_X B$$
for each commutative $k$-algebra $B$. Let us observe that
${\mathcal A}_X = \underset{X}{\prod} {\bf k}
=(\underset{X}{\oplus} {\bf k})^*$ is a commutative algebra
scheme.

Let us prove that ${\rm Spec}\, {\mathcal A}_X = {\mathcal X}$.

The morphisms of functors of $k$-algebras are, in particular,
morphisms of functors of $k$-modules $$({\rm Spec}\, {\mathcal
A}_X )(B) = {\rm Hom}_{k-alg} ({\mathcal A}_X, {\bf B} ) \subset
{\rm Hom}_k ( {\mathcal A}_X, {\bf B} )= {\rm Hom}_k ( \underset{x
\in X}{\prod} {\bf k}, {\bf B} )$$ and from \cite[4.5]{Amel} we
know that morphisms of functors of $k$-modules from
$\underset{X}{\prod} {\bf k}$ to ${\bf B}$ factorize through the
projection onto a finite number of factors. Therefore, given a
morphism $\underset{X}{\prod} {\bf k} \to {\bf B}$ there exists a
subset $Y \subset X$ of finite order, such that we get the
factorization $$\underset{X}{\prod} {\bf k} \to
\underset{Y}{\prod} {\bf k} \to {\bf B}.$$ Every continuous
morphism ${\rm Spec}\, B \to X$ has finite image, because ${\rm
Spec}\, B$ is compact and $X$ is discrete. We can assume that $X$
is a finite set.

Thus, we get

$$\aligned {\rm Hom}_{k- alg}(\prod_X {\bf k,B})& ={\rm Hom}_{k- alg}(\prod_X k,B)\\
&  = {\rm Cont. Aplic.}_{{\rm Spec}\, k}({\rm Spec}\, B, {\rm Spec}\,
\prod_X k = \coprod_X {\rm Spec}\, k) \\&  = {\rm
Aplic}_{cont.}({\rm Spec}\, B, X) . \endaligned$$

Later we will see that the ring of functions of ${\rm Spec}\,
{\mathcal A}_X$ coincides with ${\mathcal A}_X$, that is,
$${\bf Hom}_{functors}({\rm Spec}\, {\mathcal
A}_X , {\bf k}) = {\mathcal A}_X .$$
\end{eje}

\begin{pro}
Let $\mathcal A$ a functor of commutative algebras. Then, ${\rm
Spec}\,{\mathcal A}={\bf Hom}_{k-alg}({\mathcal A},{\bf k})$.
\end{pro}

\begin{proof}
By the adjoint functor formula (\cite[1.15]{Amel}) (restricted to
the morphisms of algebras) it holds that
$$\aligned {\bf Hom}_{k- alg} ({\mathcal A}, {\bf k}) (B) &
= {\rm Hom}_{B- alg} ({\mathcal A}_{|B}, {\bf B})
= {\rm Hom}_{k-alg} ({\mathcal A}, {\bf B}) \\
& = ({\rm Spec}\, {\mathcal A}) (B) . \endaligned $$
\end{proof}

Therefore, ${\rm Spec}\,{\bf A^*}={\bf Hom}_{k-alg} ({\bf A}^*,{\bf
k}) \subset {\bf Hom}_{k}({\bf A}^*,{\bf k}) = {\bf A}$.

If $E$ is a finitely-generated $k$-module, we will say that ${\bf
E}$ is coherent.

\begin{teo} \label{1.5}
Let ${\bf A}^*$ be a commutative algebra scheme and let us
assume that $A$ is a projective $k$-module. By \cite[4.12]{Amel},
${\bf A}^* = \underset{\underset{i}{\longleftarrow}}{\lim}\, {\bf
A}_i$, where ${\bf A}_i$ are the coherent algebras that are quotients
of ${\bf A}^*$. It holds that $${\rm Spec}\,{\bf A}^* =
\underset{\underset{i}{\longrightarrow}} {\lim}\, {\rm Spec}\,
{\bf A}_i.$$
\end{teo}

\begin{proof}
By \cite[4.5]{Amel}, the image of every $k$-linear morphism ${\bf
A}^* \to {\bf B}$ is a coherent $k$-module; therefore, every
morphism of functors of $k$-algebras ${\bf A}^* \to {\bf B}$
factorizes through a coherent algebra ${\bf A}_i$ that is a
quotient of ${\bf A}^*$. Then,

\begin{equation*} \begin{split} ({\rm Spec}\, {\bf A^*})(B) & ={\rm
Hom}_{k-alg} ({\bf A}^*, {\bf B}) =
\underset{\underset{i}{\longrightarrow}} {\lim}\, {\rm Hom}_{k-
alg} ({\bf A}_i, {\bf B}) \\ & =
\underset{\underset{i}{\longrightarrow}}{\lim}\,({\rm Spec}\, {\bf
A}_i) (B) .
\end{split}\end{equation*}
\end{proof}

\begin{pro}
Let ${\mathcal B}$ be a functor of commutative $k$-algebras. Then,

\begin{enumerate}

\item ${\rm Hom}_{functors}({\rm Spec}\, {\bf A},{\rm Spec}\, {\mathcal
B} ) = {\rm Hom}_{k-alg}({\mathcal B},{\bf A}).$

\item ${\rm Hom}_{functors}({\rm Spec}\, {\bf A}^*,{\rm Spec}\,
{\mathcal B} ) = {\rm Hom}_{k-alg}({\mathcal B},{\bf A}^*).$
\end{enumerate}
\end{pro}

\begin{proof}
It holds that  ${\rm Hom}_{functors} (({\rm Spec}\, A)^\cdot, F) =
F(A)$ for every functor $F$, by Yoneda's lemma.

\begin{enumerate}
\item
$$ \aligned {\rm Hom}_{functors}({\rm Spec}\, {\bf A},{\rm Spec}\, {\mathcal B})
& = {\rm Hom}_{functors} ( ({\rm Spec}\, {A})^\cdot,{\rm Spec}\, {\mathcal B}) \\
& ={\rm Hom}_{k-alg} ({\mathcal B},{\bf A}) . \endaligned$$

\item $$\aligned {\rm Hom}_{functors}({\rm Spec}\, {\bf A}^*,{\rm
Spec}\, {\mathcal B}) & ={\rm Hom}_{functors}
(\underset{\underset{i}{\longrightarrow}} {\lim}\, {\rm Spec}\, {\bf
A}_i,{\rm Spec}\, {\mathcal B}) \\ & =
\underset{\underset{i}{\longleftarrow}}{\lim}\,{\rm Hom}_{functors}
({\rm Spec}\, {\bf A}_i,{\rm Spec}\, {\mathcal B})
\\ & = \underset{\underset{i}{\longleftarrow}}
{\lim}\,{\rm Hom}_{k-alg} ({\mathcal B},{\bf A}_i) \\
& = {\rm Hom}_{k-alg} ({\mathcal B},{\bf A}^*) .
\endaligned $$\end{enumerate}
\end{proof}

\begin{pro}
For ${\mathcal A}=\bf A, A^*$, it holds that
$${\bf Hom}_{functors}({\rm Spec}\, {\mathcal A}, {\bf k}) = {\mathcal A} .$$
\end{pro}

\begin{proof}
$$ \aligned {\bf Hom}_{functors}({\rm Spec}\, {\mathcal A}, {\bf
k}) & = {\bf Hom}_{functors}({\rm Spec}\, {\mathcal A}, {\rm
Spec}\,{\bf k[x]}) \\ & = {\bf Hom}_{k-alg}({\bf k[x],
}{\mathcal A})= {\mathcal A} . \endaligned$$
\end{proof}

Explicitly, ${\mathcal A}= {\bf Hom}_{functors} ({\rm Spec}\,
{\mathcal A}, {\bf k})$, $w\mapsto \tilde{w}$, where
$\tilde{w}(\phi) = \phi(w)$, for every $\phi \in {\rm Spec}\,
{\mathcal A}= {\bf Hom}_{k-alg} ({\mathcal A}, {\bf k})$.

\section{Closure of dual modules and algebras}

\begin{defi}
We will say that a functor of $k$-modules $F$ is dual if there
exists a functor of $k$-modules $H$ such that $F \simeq H^*$.
\end{defi}

\begin{pro}
Let $F$ be a functor of $k$-modules such that $F^*$ is a reflexive
functor. The closure of dual functors of $k$-modules of $F$ is
$F^{**}$, that is, it holds the functorial equality $${\rm
Hom}_k(F,G)={\rm Hom}_k(F^{**},G)$$ for every dual functor of
$k$-modules $G$.
\end{pro}

\begin{proof}
Let us write $G=H^*$, then
$$\aligned {\rm Hom}_k(F,G) & ={\rm Hom}_k(F,H^*)={\rm Hom}_k(H,F^*)={\rm
Hom}_k(H,F^{***}) \\ & ={\rm Hom}_k(F^{**},H^*)={\rm Hom}_k(F^{**},G) . \endaligned$$
\end{proof}

\begin{defi}
Let $F$ be a functor of sets. We will denote by $k[F]$ the functor over the category of commutative $k$-algebras defined by $k[F](B)= \{$formal finite $B$-linear combinations of $f \in F(B) \}$.
\end{defi}

\begin{cor}\label{2.3}
Let $\mathcal A$ be a functor of commutative $k$-algebras. Let us
assume that ${\bf Hom}_{functors}({\rm Spec}\, {\mathcal A}, {\bf
k})= \mathcal A$ and that $\mathcal A$ is a reflexive functor of
$k$-modules.Then, the closure of dual functors of $k$-modules of
$k[{\rm Spec}\, {\mathcal A}]$ is ${\mathcal A}^*$, that is,
$${\rm Hom}_{functors}({\rm Spec}\, {\mathcal A}, F) = {\rm Hom}_k(k[{\rm
Spec}\, {\mathcal A}], F) = {\rm Hom}_k({\mathcal A}^*, F)$$ for
every dual functor of $k$-modules $F$.
\end{cor}

\begin{proof}
As $k[{\rm Spec}\, {\mathcal A}]^*={\bf Hom}_{functors}({\rm
Spec}\, {\mathcal A}, {\bf k})=\mathcal A$, the closure of dual
functors of $k$-modules of $k[{\rm Spec}\, {\mathcal A}]$ is
$k[{\rm Spec}\, {\mathcal A}]^{**}={\mathcal A}^*$.
\end{proof}

\begin{pro}\label{2.4}
Let $F$ be a functor of $k$-algebras such that $F^*$ is a
reflexive functor of $k$-modules. The closure of dual functors of
$k$-algebras of $F$ is $F^{**}$, that is, it holds the functorial
equality $${\rm Hom}_{k-alg}(F,G)={\rm Hom}_{k-alg}(F^{**},G)$$
for every dual functor of $k$-algebras $G$. In particular, ${\rm
Spec}\, F= {\rm Spec}\, F^{**}$.
\end{pro}

\begin{proof}
Let us observe that
$$\aligned (F \otimes \overset{n}{\ldots} \otimes F)^*  =
& {\rm Hom}_k(F \otimes \overset{n}{\ldots} \otimes F,{\bf k})
= {\rm Hom}_k (F,(F \otimes \overset{n-1}{\ldots} \otimes F)^*)\\
\underset{Induction}= & {\rm Hom}_k (F, (F^{**} \otimes
\overset{n-1}{\ldots} \otimes F^{**})^*) = {\rm Hom}_k
(F^{**},(F^{**} \otimes \overset{n-1}{\ldots} \otimes F^{**})^*)\\
= & (F^{**} \otimes \overset{n}{\ldots} \otimes F^{**})^*
.\endaligned$$ Therefore, given a dual functor of $k$-modules
$S^*$, $$\aligned {\rm Hom}_k (F \otimes \ldots \otimes F, S^*) &
={\rm Hom}_k(S, (F \otimes \ldots \otimes F)^*) \\ & = {\rm
Hom}_k(S, (F^{**} \otimes \ldots \otimes F^{**})^*)\\ & = {\rm
Hom}_k (F^{**} \otimes \ldots \otimes F^{**}, S^*) .\endaligned$$

If we consider $S^*=F^{**}$, it follows easily that the structure of
algebra of $F$ defines a structure  of algebra on $F^{**}$. Finally,
if we consider $S^*= G$, we obtain that ${\rm Hom}_{k-alg}(F,G)={\rm
Hom}_{k-alg}(F^{**},G)$.
\end{proof}

\begin{com}
Let us observe that if $F$ is also a functor of commutative
algebras, then so is $F^{**}$: since ${\rm Hom}_k (F \otimes F,
F^{**}) = {\rm Hom}_k (F^{**} \otimes F^{**}, F^{**})$, the
morphism $F \otimes F \to F$, $f \otimes f'\mapsto ff'-f'f=0$
extends to a unique morphism $F^{**} \otimes F^{**} \to F^{**}$
(which is $f \otimes f'\mapsto ff'-f'f =0 $).
\end{com}

Let ${\bf A}^*$ and ${\bf B}^*$ be commutative $k$-algebra
schemes. By \cite[1.8]{Amel}, we have that $({\bf A}^* \otimes_k
{\bf B^*})^* = {\bf Hom}_k ({\bf A^*} \otimes_k {\bf B^*},{\bf k})
= {\bf Hom}_k ({\bf A^*},{\bf B}) = {\bf A} \otimes_k {\bf B}$.
Then, ${\rm Spec}\,{\bf A}^* \times {\rm Spec}\,{\bf B}^*= {\rm
Spec}\, ({\bf A}^* \otimes_k {\bf B}^*)
\stackrel{\text{\ref{2.4}}}{=} {\rm Spec}\, ({\bf A} \otimes_k
{\bf B})^* $.

Let $\mathcal A$ be a functor of $k$-algebras and let us assume
that $G= {\rm Spec}\, \mathcal A$ is a functor of monoids. The
functor of $k$-modules $k[G]$ is obviously a functor of
$k$-algebras. Given a functor of $k$-algebras $\mathcal B$, it is
easy to check the equality $${\rm Hom}_{monoids} (G, {\mathcal B})
= {\rm Hom}_{k-alg} (k[G], \mathcal B) .$$

\begin{pro}\label{2.5}
Let $G= {\rm Spec}\, \mathcal A$ be a functor of monoids. Let us
assume that ${\bf Hom}_{functors}(G , {\bf k}) = \mathcal A$ and
that $\mathcal A$ is a reflexive functor of $k$-modules. Then the
closure of dual functors of algebras of $k[G]$ coincides with
$\mathcal A^*$, that is,

$${\rm Hom}_{monoids}(G, {\mathcal B}) = {\rm Hom}_{k-alg}(k[G],
{\mathcal B}) = {\rm Hom}_{k-alg}({\mathcal A^*}, {\mathcal B})$$
for every dual functor of $k$-algebras $\mathcal B$.
\end{pro}

\begin{proof}
$k[G]^*=\mathcal A$ is reflexive, then the closure of dual functors of
algebras of $k[G]$ is $\mathcal A^*$.
\end{proof}

\begin{teo}\label{Gmodulos}
Let $G= {\rm Spec}\, \mathcal A$ be a functor of monoids. Let us
assume that ${\bf Hom}_{functors} ({\rm Spec}\, {\mathcal A}, {\bf
k}) = \mathcal A$ and that $\mathcal A$ is a reflexive functor of
$k$-modules. The category of $G$-modules is equi\-valent to the
category of ${\mathcal A}^*$-modules. Likewise, the category of
dual functors of $G$-modules is equivalent to the category of dual
functors of ${\mathcal A}^*$-modules.
\end{teo}

\begin{proof}
Let $E$ be a $k$-module. Let us observe that ${\bf End}_k ({\bf
E}) = ({\bf E}^* \otimes {\bf E})^*$ is a dual functor. Therefore,
$${\rm Hom}_{k-alg} (k[G],{\bf End}_k({\bf E})) = {\rm
Hom}_{k-alg} ({\mathcal A}^*, {\bf End}_k({\bf E})) . $$ In
conclusion, endowing $E$ with a structure of $G$-module is
equivalent to endowing $E$ with a structure of ${\mathcal
A}^*$-module.

Let us observe that ${\rm Hom}_k (k[G], {\bf E}) = {\rm Hom}_k
({\bf E}^*, {\mathcal A})= {\rm Hom}_k ({\mathcal A}^*, {\bf E})$.
Hence, given any two $G$-modules (or ${\mathcal A}^*$-modules)
$E$, $E'$, a linear morphism $f: E \to E'$ and $e \in E$, we will
have that the morphism $f_1: k[G] \to {\bf E}'$, $f_1(g) := f(ge)
- gf(e)$ is null if and only if the morphism $f_2: {\mathcal A}^*
\to {\bf E}'$, $f_2(a) := f(ae) - a f(e)$ is null. In conclusion,
${\rm Hom}_{G-mod} ({\bf E}, {\bf E'}) = {\rm Hom}_{{\mathcal
A}^*} ({\bf E}, {\bf E'})$.
\end{proof}

Let us highlight that the structure of functor of algebras of
$\mathcal A^*$ is the one that makes the embedding ${\rm
Spec}\,{\mathcal A}\hookrightarrow \mathcal A^*$ be a morphism of
monoids.

\begin{eje}
The ${\mathbb C}$-linear representations of $\Z$ are equivalent to
the ${\mathbb C}[\Z]$-modu\-les. ${\mathbb C}[\Z]={\mathbb
C}[x,1/x]$, $n\mapsto x^n$, is a principal ideal domain. Thus, if
$E$ is a finite ${\mathbb C}$-linear representation of $\Z$, then
$$E= \underset{\alpha\neq 0, n_\alpha>0}{\oplus} {\mathbb
C}[x]/(x-\alpha)^{n_\alpha}$$ such that $n\cdot (\bar
p_{n_\alpha})_{n_\alpha}=(\overline{x^n\cdot
p_{n_\alpha}})_{n_\alpha}$.
\end{eje}

\begin{pro}\label{2.9}
Let $G = {\rm Spec}\, {\mathcal A}$ be a functor of monoids. Let
us assume that $ {\bf Hom}_{functors} ( {\rm Spec}\, {\mathcal A},
{\bf k}) = {\mathcal A}$ and that ${\mathcal A}$ is a reflexive
functor of $k$-modules. Let ${\mathcal V}$ be a dual functor of
$G$-modules (on the left and on the right). Then, $$Der_k
({\mathcal A}^*, {\mathcal V}) = Der(G,{\mathcal V})$$ where
$Der(G, {\mathcal V}) := \{ D \in {\rm Hom}_{functors} (G, {\mathcal
V}) : D(g \cdot g') = D(g) \cdot g' + g \cdot D(g') \}$.
\end{pro}

\begin{proof}
$$ \aligned Der_k ({\mathcal A}^*, {\mathcal V}) & =
\{ f \in {\rm Hom}_{k-alg} ({\mathcal A}^*, {\mathcal A}^*
\oplus {\mathcal V} \epsilon), \epsilon^2 = 0, f(a)= a \, {\rm mod}
(\epsilon) \} \\ & = \{ f \in  {\rm Hom}_{monoids} (G, {\mathcal A}^*
\oplus {\mathcal V} \epsilon), \epsilon^2 =0,  f(g) =
g \, {\rm mod} (\epsilon) \} \\ & = Der(G, {\mathcal V}) . \endaligned$$
\end{proof}

\section{Functorial Cartier Duality}

\begin{defi}
Given a functor of abelian monoids $G$, where we regard ${\bf k}$
as a monoid with the multiplication operation,  $$G^* := {\bf
Hom}_{monoids} (G,{\bf k}) $$ will be referred to as the dual
monoid of $G$.
\end{defi}

If $G$ is a functor of groups, then $G^* = {\bf Hom}_{groups} (G, G_m^{\cdot})$.

\begin{teo}\label{dual}
Assume that ${\rm Spec}\, {\mathcal  A}$ is a functor of abelian
monoids, $\mathcal A$ is a reflexive functor of $k$-modules and ${\bf
Hom}_{functors}({\rm Spec}\, {\mathcal  A},{\bf k})=\mathcal A$. Then

$${\rm Spec}\, {\mathcal  A}^*=({\rm Spec}\, {\mathcal  A})^*$$
(in particular, this equality shows that ${\rm Spec}\,\mathcal
A^*$ is a functor of abelian monoids).
\end{teo}

\begin{proof}
$$\aligned {\rm Spec}\, {\mathcal  A}^* & = {\bf Hom}_{k-alg}
({\mathcal  A}^*,{\bf k}) = {\bf Hom}_{monoids} ({\rm Spec}\,
{\mathcal A},{\bf k})\\ & = ({\rm Spec}\, {\mathcal A})^* .\endaligned $$
\end{proof}

\begin{com}\label{comentario}
Explicitly, ${\rm Spec}\, {\mathcal A}^* = {\bf Hom}_{monoids} ({\rm
Spec}\, {\mathcal A},G_m^\cdot)$, $\phi \mapsto \tilde{\phi}$,
where $\tilde{\phi} (x)$ $ = \phi(x)$, for every $\phi \in {\rm
Spec}\, {\mathcal A}^* = {\bf Hom}_{k-alg} ({\mathcal  A}^*,{\bf
k})$ and $x \in {\rm Spec}\, {\mathcal A} \subset \mathcal A^*$.

Let $G= {\rm Spec}\, {\mathcal A}$. $G^*$ is a functor of abelian
monoids ($(f \cdot f')(g) := f(g) \cdot f'(g)$, for every $f,f'\in
G^*$ and $g \in G$), the inclusion $G^* = {\bf Hom}_{monoids}
(G,{\bf k}) \subset {\bf Hom}_{functors} (G,{\bf k})={\mathcal A}$
is a morphism of monoids and the diagram $$\xymatrix{ ({\rm
Spec}\, {\mathcal A})^* \ar@{^{(}->}[r] \ar@{=}[d] & {\bf
Hom}_{functors}({\rm Spec}\,{\mathcal A}, {\bf k}) = \mathcal A
\ar@{=}[d] \\ {\rm Spec}\, {\mathcal A^*} \ar@{^{(}->}[r] & {\bf
Hom}_k({\mathcal A}^*,{\bf k})=\mathcal A^{**}}$$ is commutative.
\end{com}

\begin{teo}\label{Cartier}
Assume that $G={\rm Spec}\, {\mathcal  A}$ is a functor of abelian
monoids, $\mathcal A$ is a reflexive functor of $k$-modules and ${\bf
Hom}_{functors}({\rm Spec}\, {\mathcal  A},{\bf k})=\mathcal A$, ${\bf
Hom}_{functors}({\rm Spec}\, $ ${\mathcal  A}^*, {\bf k}) = \mathcal A^*$. Then:
\begin{enumerate}

\item The morphism $G\overset{**}\to G^{**}$, $g \mapsto
g^{**}$, where $g^{**} (f) := f(g)$ for every $f \in G^* $, is an
isomorphism.

\item ${\rm Hom}_{monoids} (G_1,G_2) = {\rm Hom}_{monoids}
(G_2^*,G_1^*)$.
\end{enumerate}
\end{teo}

\begin{proof} $ $
\begin{enumerate}
\item It is easy to check that the diagram
$$\xymatrix{{\rm Spec}\, {\mathcal A^{**}} \ar[r]^-\sim_-{\text{\ref{dual}}} & ({\rm
Spec} \, {\mathcal A^*})^*
& ({\rm Spec}\, {\mathcal A})^{**}
\ar[l]_-\sim^-{\text{\ref{dual}}}\\ & {\rm Spec}\, {\mathcal A}
\ar@{=}[ul] \ar[ur]^-{**}& }$$ is commutative. Hence the morphism $**$ is an isomorphism.

\item Every morphism of monoids $G_1 \to G_2$, taking ${\bf
Hom}_{monoids} (-,{\bf k})$, defines a morphism $G_2^* \to G_1^*$.
Taking ${\bf Hom}_{monoids}(-,{\bf k})$ we get the original
morphism $G_1 \to G_2$, as it is easy to check. Likewise, every
morphism of monoids $G_2^* \to G_1^*$, taking ${\bf
Hom}_{monoids}(-,{\bf k})$, defines a morphism $G_1 \to G_2$.
Taking ${\bf Hom}_{monoids}(-,{\bf k})$ we get the original
morphism $G_2^* \to G_1^*$.
\end{enumerate}
\end{proof}

\begin{teo}
The category of abelian affine $k$-monoids $G={\rm Spec}\, A$ is
anti-equiva\-lent to the category of functors ${\rm Spec}\, {\bf
A}^*$ of abelian monoids (we assume the $k$-modules $A$ are
projective).
\end{teo}

In particular, we get the Cartier duality for finite commutative
algebraic groups (\cite[\S 9.9]{carlos}). Dieudonn\'{e}
(\cite[Ch. I, \S 2, 13]{dieudonne}) proves the equivalence between
the category of $k$-coalgebras in groups and the dual category of
linearly compact $k$-algebras in cogroups (where $k$ is a field).

\begin{pro}
Let $G_1$ and $G_2$ be a pair of functors of abelian monoids. It holds
that $$(G_1 \times G_2)^*= G_1^* \times G_2^* .$$
\end{pro}

\begin{proof}
$$\aligned (G_1 \times G_2)^* & ={\bf Hom}_{monoids} (G_1\times G_2, {\bf k})
\\ & = {\bf Hom}_{monoids} (G_1, {\bf k}) \times {\bf Hom}_{monoids} (G_2, {\bf k})
= G_1^* \times G_2^*. \endaligned$$
\end{proof}

Let us show some examples.

\begin{eje}
Let $k$ be a field. Let $\Z := {\rm Spec}\, \underset{\Z}{\prod} {\bf k}$.  Obviously $\Z^* = G_m$, therefore $G_m^* = \Z$. In other words, ${\bf Hom}_{grupos} (G_m,G_m)=\Z$: given $\tau : G_m \to G_m$, there
exists $n \in \Z$ such that $\tau(\alpha) = \alpha^n$.

Given a quasi-coherent $\prod_X {\bf k}$-module ${\bf E}$, then ${\bf E} = \oplus_{x \in X} {\bf E}_x$ in the obvious way. Given a morphism of algebra schemes $\prod_X {\bf k} \to {\bf A}^*$, then ${\bf A}^* = \prod_{x \in X} {\bf A}^*_x$. Moreover, if the morphism ${\rm Spec}\, {\bf A}^* \to {\rm Spec}\, \prod_X {\bf k}$ is injective, then ${\bf A}^*_x = {\bf k}$ or ${\bf A}^*_x =0$ for each $x \in X$. Therefore, ${\rm Spec}\, {\bf A}^* \subset {\rm Spec}\, \underset{\Z}{\prod} {\bf k} = \Z$ is a subgroup if and only if ${\rm Spec}\, {\bf A}^* =n \cdot \Z$. Dually, the quotient algebraic groups of $G_m$ are the epimorphisms $G_m \to G_m$, $t \mapsto t^n$. Hence, the subgroups of $G_m$ are the $\mu_n = \{
\alpha \in G_m : \alpha^n = 1 \} = {\rm Spec}\, k[x]/(x^n-1)=(\Z/(n))^*$.

Likewise, $(G_m \times \overset{n}{\ldots} \times G_m)^* = \Z
\times \overset{n}{\ldots} \times \Z$. By the theory of abelian
groups (or $\Z$-modules), the subgroups of $\Z \times
\overset{n}{\ldots} \times \Z$ are isomorphic to $\Z^r$ through an
injective morphism $\phi$ (let us say that its matrix is
$(n_{ij})$), $$\Z \times \overset{r}{\ldots} \times \Z
\overset{\phi}{\to} \Z \times \overset{n}{\ldots} \times \Z$$
whose cokernel is isomorphic to $\Z^{n-r} \times \Z/(n_1) \times
\ldots \times \Z/(n_r)$. Dually, the subgroups of $G_m \times
\overset{n}{\ldots} \times G_m$ are isomorphic to $G_m^{n-r}
\times \mu_{n_1} \times \ldots \times \mu_{n_r}$, that is
isomorphic to the kernel of the epimorphism $$G_m^n \to G_m^r,\;
(t_1, \ldots ,t_n) \mapsto (t_1^{n_{11}} \cdot \ldots \cdot
t_n^{n_{n1}},\ldots, t_1^{n_{1r}} \cdot \ldots \cdot t_n^{n_{nr}})
. $$
\end{eje}

\begin{pro}\cite[Ch. III, 8.12]{borel}
The category of diagonalizable algebraic groups is anti-equivalent
to the category of $\Z$-modules of finite type (or
finitely-generated abelian groups).
\end{pro}

\begin{eje}[{\bf Affine toric varieties}]
Let $M$ be a set with structure of abelian (multiplicative)
monoid. The constant functor $\mathcal{M} = {\rm Spec}\, \prod_M
{\bf k}$ is a functor of abelian monoids. The dual functor is the
abelian $k$-monoid ${\mathcal M}^* = {\rm Spec}\, \oplus_M {\bf k}
= {\rm Spec}\, k[M]$.

We will say that an abelian monoid $M$ is classic if it is finitely generated and its associated group $G$ is torsion-free and the natural morphism $M \to G$ is injective. It is easy to prove that $M$ is classic if and only if $k[M]$ is an integral finitely generated $k$-algebra.
\\ {\bf Theorem. } The category of abelian monoids (respectively finitely generated, classic) is anti-equivalent to the category of affine $k$-schemes of semisimple abelian monoids (respectively algebraic, integral algebraic).

As $G = \mathbb{Z}^n$, the morphism $M \to G$ induces a morphism $G^n_m \to {\mathcal M}^*$. In particular, $G^n_m$ operates on ${\mathcal M}^*$. Furthermore, as $k[G]$ is the localization of $k[M]$ by the algebraically closed system $M$, the morphism $G^n_m \to {\mathcal M}^*$ is an open injection. We will say that an integral affine algebraic variety on which the torus operates with a dense orbit is an affine toric variety. It is easy to prove that there exists a one-to-one correspondence between affine toric varieties with a fixed point whose orbit is transitive and dense, and classic monoids.
\end{eje}

\begin{eje} Let $G_a={\rm Spec}\, k[x]$ be the additive group.
Let us compute $G_a^*$. Assume ${\rm car}\, k=0$. Let $\bf k[[z]]$
be the $k$-algebra scheme defined by ${\bf k[[z]]}(B):=B[[z]]$. Let $\{ \omega_0, \ldots, \omega_n, \ldots \} \subset k[x]^*$ be the dual basis of $\{ 1, x, \ldots, x^n , \ldots \} \subset k[x]$.
Let $\phi\colon \bf k[x]^*\to k[[z]]$ be the linear isomorphism
defined by 
$\phi(\omega_i) = \dfrac{z^i}{i!}$.
The composition of the natural morphism $G_a\to \bf k[x]^*$ with
$\phi$ is equal to the morphism $G_a\to \bf k[[z]]$,
$\alpha\mapsto e^{\alpha z}$, which is a morphism of monoids.
Therefore, $\phi$ is an isomorphism of functors of algebras.

$$\aligned G_a^*(B) & = {\rm Hom}_{k-alg} ({\bf k[x]}^*,
{\bf B}) = {\rm Hom}_{k-alg} ({\bf k[[z]]}, {\bf B}) \\
& = \underset{\underset{n} {\longrightarrow}}{\lim}\, {\rm
Hom}_{k-alg} ({k[z]/(z^n)},B) = {\rm rad}\, B .
\endaligned$$ Let us follow the notation $G_a^*={\rm rad}\, {\bf k}$.
Explicitly, we assign to $n \in {\rm rad}\, {\bf k}$ the morphism
$G_a \to G_m$, $\alpha \mapsto e^{\alpha \cdot n}$. Finally,
$({\rm rad}\, {\bf k})^*=G_a$; explicitly, $\alpha \in G_a$ defines the morphism ${\rm rad}\, {\bf k} \to G_m$, $u \mapsto e^{\alpha \cdot u}$.
\end{eje}

\begin{pro}
Let $G = {\rm Spec}\, A$ be an abelian $k$-monoid. Let $p \in
G^*(k)$ and let us regard $k$ as a $G$-module via $p$. Then $$ T_p
G^* := Der_k ({\bf A}^*, {\bf k})= Der(G, {\bf k}) .$$

Let $q \in G(k)$ and let us regard $k$ as a $G^*$-module via $q$.
Then $$T_q G = Der_k ({\bf A}, {\bf k}) = Der(G^*, {\bf k}) .$$
\end{pro}

\begin{proof}
It is a consequence from Proposition \ref{2.9}, where ${\mathcal V}
= {\bf k}$.
\end{proof}

Let ${\bf I}_p^* := \ker p$, then ${\rm Hom}_k ( {\bf I}_p^*/{\bf
I}^{*2}_p , {\bf k}) = Der_k ({\bf A}^*, {\bf k}) = 
Der (G, {\bf k})$.

\section{Structure of commutative $k$-groups in characteristic zero}

For simplicity we assume that $k$ is an algebraically closed
field. Let ${\bf A}^*$ be a $k$-scheme of commutative algebras and ${\bf
M}^* = \prod_J {\bf k}$ the maximal semisimple quotient of ${\bf
A}^*$. The maximal spectrum of a scheme of algebras is defined to
be the set of its maximal bilateral ideal schemes and it is the
same as the set of isomorphism classes of simple ${\bf
A}^*$-modules (see \cite[6.7]{Amel}). Then, ${\rm Spec}_{max} {\bf
A}^* = {\rm Spec}_{max} {\bf M}^* = J$. ${\bf A}^*$ is the inverse
limit of its cokernels of finite dimension. As every finite
commutative $k$-algebra is a direct product of local $k$-algebras,
it is easy to prove that ${\bf A}^* = \prod_{j \in J} {\bf
A}^*_j$, where the maximal semisimple quotient of ${\bf A}^*_i$ is
${\bf k}$. Then, there exists an only section of algebra
$k$-schemes $s: {\bf M}^* \to {\bf A}^*$ of the epimorphism ${\bf
A}^* \to {\bf M}^*$.

The maximal semisimple quotient of ${\bf A}^* \bar{\otimes} {\bf
A}^* = ( {\bf A} \otimes {\bf A})^*$ is ${\bf M}^* \bar{\otimes}
{\bf M}^* = ( {\bf M} \otimes {\bf M})^*$ (\cite[6.23]{Amel}). If
${\rm Spec}\, {\bf A}^*$ is also a functor of abelian groups, then
${\rm Spec}\, {\bf M}^*$, which is the constant functor $J$, is a
functor of abelian subgroups. Furthermore, the morphism $s^*: {\rm
Spec}\, {\bf A}^* \to {\rm Spec}\, {\bf M}^*$ induced by $s$ is a
morphism of functors of groups. Therefore, if $1 \in J$ is the
identity element, we have the following decomposition of group
functors $$ {\rm Spec}\, {\bf A}^* = {\rm Spec}\, {\bf M}^* \times
{\rm Spec}\, {\bf A}_1^*.$$ In particular, $G = {\rm Spec}\, A$,
which is the dual of $G^* = {\rm Spec}\, {\bf A}^*$, is a direct
product of a diagonalizable group (every linear representation is
a direct sum of simple modules of rank $1$) and a unipotent group
(every linear representation has got a filtration of invariant
factors).

Given a $k$-vector space scheme ${\bf E}^*$, let us denote by
$\bar{S}^n {\bf E}^*$ the representant of the functor $$F({\bf
V}^*)= {\rm Sim}_k ({\bf E}^* \times \stackrel{n}{\ldots} \times
{\bf E}^*, {\bf V}^* ) .$$ $\bar{S}^n {\bf E}^*$ is the quotient
of $ {\bf E}^* \bar{\otimes} \stackrel{n}{\ldots} \bar{\otimes}
{\bf E}^* = ({\bf E} \otimes \stackrel{n}{\ldots} \otimes {\bf
E})^*$ by the minimum $k$-vector space subscheme that contains
elements of the following type $$\ldots \otimes w \otimes \ldots
\otimes w' \otimes \ldots - \ldots \otimes w' \otimes \ldots
\otimes w \otimes \ldots \; .$$ $\bar{S}^n {\bf E}^*$ coincides
with the closure of $k$-vector space  schemes of $S^n {\bf E}^*$
and it coincides with $(( {\bf E} \otimes \stackrel{n}{\ldots}
\otimes {\bf E})^{S_n})^*$. In characteristic zero, $(E \otimes
\stackrel{n}{\ldots} \otimes E)^{S_n} = S^n E$, then $\bar{S}^n
{\bf E}^* = (S^n {\bf E})^*$.

If $E = \oplus_X k$, then ${\bf E}^* = \prod_X {\bf k}$ and
$\bar{S}^n {\bf E}^* = \prod_{S^n X} {\bf k}$.

We will denote $k[[{\bf E}^*]] := \prod_{n \in \N} \bar{S}^n {\bf
E}^*$. If ${\rm car}\, k=0$ and $G = {\bf V}^*$, then $G^* = {\rm
Spec}\, k[[{\bf V}^*]]$. Specifically, the natural morphism ${\bf
V}^* \to k[[{\bf V}^*]]$ assigns $w \in {\bf V}^*$ to $e^w$.

\begin{defi}
The kernel of the morphism ${\bf A}^* \to {\bf M}^*$ (the morphism of a scheme of algebras onto its maximal semisimple quotient) is said to be the radical of ${\bf A}^*$. We will say that ${\bf A}^*$ is local if ${\bf M}^* = {\bf k}$.
\end{defi}

\begin{lema}
Let $f: {\bf A}^* \to {\bf B}^*$ be a morphism between local
$k$-algebra schemes and let ${\bf I}^*_A$ and ${\bf I}^*_B$ be the
radical ideals of ${\bf A}^*$ and ${\bf B}^*$, respectively. If
the induced morphism ${\bf I}^*_A/{\bf I}^{*2}_A \to {\bf
I}^*_B/{\bf I}^{*2}_B$ is surjective, then $f$ is a surjective
morphism. If ${\bf B}^* = k[[{\bf E}^*]]$ and the morphism ${\bf
I}^*_A/{\bf I}^{*2}_A \to {\bf I}^*_B/{\bf I}^{*2}_B$ is
injective, then $f$ is an injective morphism.
\end{lema}

\begin{proof}
By abuse of notation, we have regarded ${\bf I}^{*n}$ as the closure of
vector space schemes of ${\bf I}^{*n}$ in ${\bf A}^*$ (as we did
in \cite{Amel}). It is convenient to denote by $({\bf I}^{*n})'$
this closure in this paragraph. The natural morphism $S^n
{\bf I}^*/{\bf I}^{*2} \to {\bf I}^{*n}/{\bf I}^{*n+1}$ is surjective
and the closure of $k$-vector space schemes of the image of ${\bf
I}^{*n}/{\bf I}^{*n+1}$ in ${\bf A}^*/({\bf I}^{*n+1})'$ is $({\bf
I}^{*n})' / ({\bf I}^{*n+1})'$. Therefore, the morphism $\bar{S}^n
({\bf I}^*/{\bf I}^{*2}) \to ({\bf I}^{*n})' / ({\bf I}^{*n+1})'$
is surjective.

If ${\bf I}^*_A/{\bf I}^{*2}_A \to {\bf I}^*_B/{\bf I}^{*2}_B$ is
a surjective morphism, then the morphisms $\bar{S}^n ({\bf
I}^*_A/{\bf I}^{*2}_A) \to \bar{S}^n ({\bf I}^*_B/{\bf I}^{*2}_B)$
are surjective and from the commutative diagram
$$\xymatrix{ \bar{S}^n ({\bf I}^*_A/{\bf I}^{*2}_A) \ar[r]^-{epi} \ar[d]^-{epi}
& \bar{S}^n ({\bf I}^*_B/{\bf I}^{*2}_B) \ar[d]^-{epi}
\\ {\bf I}^{*n}_A/{\bf I}^{*n+1}_A \ar[r] &
{\bf I}^{*n}_B/{\bf I}^{*n+1}_B }$$ it is deduced that the
morphisms ${\bf I}^{*n}_A/{\bf I}^{*n+1}_A \to {\bf I}^{*n}_B/{\bf
I}^{*n+1}_B$ are surjective. Then the morphisms ${\bf A}^*/{\bf
I}^{*n+1}_A \to {\bf B}^*/{\bf I}^{*n+1}_B$ are surjective. Taking
inverse limits (which is the dual of a direct limit of injections)
we get the epimorphism $f: {\bf A}^* \to {\bf B}^*$ (\cite[6.17 (2)]{Amel}).

If ${\bf I}^*_A/{\bf I}^{*2}_A \to {\bf I}^*_B/{\bf I}^{*2}_B$ is
an injective morphism, then the morphisms $\bar{S}^n ({\bf
I}^*_A/{\bf I}^{*2}_A) \to \bar{S}^n ({\bf I}^*_B/{\bf I}^{*2}_B)$
are injective and from the commutative diagram $$\xymatrix{
\bar{S}^n ({\bf I}^*_A/{\bf I}^{*2}_A) \ar@{^{(}->}[r]
\ar[d]^-{epi} & \bar{S}^n ({\bf I}^*_B/{\bf I}^{*2}_B) \ar@{=}[d]
\\ {\bf I}^{*n}_A/{\bf I}^{*n+1}_A \ar[r] & {\bf I}^{*n}_B/{\bf
I}^{*n+1}_B }$$ it is deduced that the morphisms ${\bf
I}^{*n}_A/{\bf I}^{*n+1}_A \to {\bf I}^{*n}_B/{\bf I}^{*n+1}_B$
are injective. Let $0 \neq a \in {\bf A}^*$ and let $n$ be such
that $a \in {\bf I}_A^{*n}$ and $a \neq {\bf I}_A^{*n+1}$ (see \cite[6.17 (2)]{Amel}). Then $f(a)\neq 0$, since $0 \neq \bar{a} \in {\bf I}^{*n}_A/{\bf I}^{*n+1}_A$ and $0 \neq \overline{f(a)} \in {\bf I}^{*n}_B/{\bf I}^{*n+1}_B$.
\end{proof}

\begin{teo}
Let $G= {\rm Spec}\, A$ be an abelian $k$-group and let $k$ be an
algebraically closed field of characterisic zero. Then, there
exists an isomorphism of $k$-groups $$G = (\underset{X}{\prod}
G_a) \times D$$ where $D$ is a diagonalizable group.
\end{teo}

\begin{proof}
We must prove that if $G$ is unipotent, then it is isomorphic to a
direct product of additive groups.

Thus, we have that ${\bf A}^*$ is local. Let ${\bf I}^*$ be its
radical. The morphism $i: G \to 1 + {\bf I}^*/{\bf I}^{*2}$,
$i(g) = 1 + \overline{g-1}$ is a morphism of $k$-groups. Let us
prove that it is an isomorphism. Consider ${\bf k}$ as the trivial $G$-module. Then, $E := Der(G, {\bf k}) = {\rm Hom}_{lin} (G, {\bf
k})$. The ring of functions of $1 + {\bf I}^*/{\bf I}^{*2} = {\bf
I}^*/{\bf I}^{*2}$ is $\oplus_n S^n E$. The morphism between the
rings induced by $i$ is the natural morphism $\oplus_n S^n E \to
A$. Dually, we have the morphism ${\bf A}^* \to k[[{\bf E}^*]] =:
{\bf B}^*$. The induced morphism ${\bf I}^*/{\bf I}^{*2} \to {\bf
I}_B^*/{\bf I}_B^{*2} = {\bf E}^*$ is the identity morphism. By
the previous lemma, ${\bf A}^* = k[[{\bf E}^*]] $ and $G= {\bf
I}^*/{\bf I}^{*2}$.
\end{proof}

\section{Inversion formula}

Let ${\bf B} \subseteq \prod_J {\bf k}$ be a quasi-coherent $k$-algebra (with or without unit). We define ${\rm Spec}\, {\bf B}$ to be the functor over the category of commutative algebras (with unit) $$({\rm Spec}\, {\bf B})(C) := \{ f \in {\rm Hom}_{k-alg} (B,C)\, \colon\, ({\rm Im} f )= C \} .$$

If $k$ is a field it is easily seen that ${\rm Spec}\, {\bf B} = {\rm Spec}\, \prod_J {\bf k}$ if and only if ${\bf B} = \oplus_J {\bf k}$.

\begin{defi}
We will say that $\oplus_J k $ is the cuasi-coherent ring of
functions of ${\rm Spec}\, \prod_J {\bf k}$.
\end{defi}

As well, $\oplus_J {\bf k}$ is characterized by being the maximal quasi-coherent $\prod_J {\bf k}$-sub\-mod\-ule of $\prod_J {\bf k}$.

Let $G^* = {\rm Spec}\, \prod_J {\bf k}$ be a functor of discrete
groups. $G^*$ operates on $G^*$ by translations on the left, then
it operates on $\prod_J {\bf k}$, then on $\oplus_J k$. The
mapping $w_{G^*}: \oplus_J k \to k$, $w_{G^*}((\lambda_i)) =
\sum_i \lambda_i$ is the only $G^*$-invariant linear form up to a
multiplicative factor. We will say that the linear mapping
$$\phi_{G^*}: \oplus_J {\bf k} \to (\prod_J {\bf k})^*
\hookrightarrow ( \oplus_J {\bf k} )^*, a \mapsto w_{G^*}(a \cdot
-)$$ is the transposed Fourier transform of $G^*$.

$G = {\rm Spec}\, A$ is a semisimple $k$-group if and only if
there exists a (unique) linear mapping $w_G: A \to k$,
$G$-invariant and such that $w_G(1)=1$ (see \cite[2.10]{Amel2}).
The morphism $\phi_G: A \to A^*$, $\phi_G(a) := w_G(a^* \cdot -)$
(where $*: A \to A$ is the morphism between rings induced by the
morphism $G \to G, g \mapsto g^{-1}$) is called the Fourier
transform of $G$ (see \cite{Amel3}). If $G$ is also commutative,
then $A^* = \prod_J k$ and ${\rm Im}\, \phi_G = \oplus_J k =:A'$.
Let $e: A \to k$ the identity element of $G$ (that is, the unit of
$A^*$). The diagram $$\xymatrix{ A \ar[dr]^-{e}
\ar@{=}[rr]^-{\phi_G} & & \ar[dl]^-{w_{G^*}}  A' \\ & k & }$$ is commutative (see \cite{Amel3}, section 2).

\begin{teo}
The transposed Fourier transform of the dual group of a semisimple abelian $k$-group is the inverse of the Fourier transform of the
semisimple abelian group.
\end{teo}

\begin{proof}
$\phi_{G^*}(\phi_G(a))(w) = w_{G^*}(\phi_G(a) \cdot w) =
w_{G^*}(\phi_G(a \cdot w)) = (a \cdot w)(e) = a(w)$ for all $a \in
A$ and $w \in A^*$.
\end{proof}

\end{document}